\documentclass[12pt,a4paper]{amsart}
\usepackage{amssymb,latexsym}

\def\N{\mathbb N}
\def\Z{\mathbb Z}
\def\R{\mathbb R}
\def\D{\mathbb D}

\def\T{\mathbb T}

\def\H2{{\mathcal H}^2(\D)}
\newtheorem{lem}{LEMMA}
\newtheorem{theo}[lem]{THEOREM}

\newtheorem{coro}[lem]{COROLLARY}
\newtheorem{prop}[lem]{PROPOSITION}
\newtheorem{proposition}[lem]{PROPOSITION}

\newtheorem{remark}[lem]{Remark}

\begin{document}
\title[ Multilinear Hankel operators]{Truncations of Multilinear Hankel operators}
\author{Aline Bonami, Sandrine Grellier \& Mohammad Kacim}
\address{MAPMO\\Universit\'e d'Orl\'eans\\ Facult\'e des Sciences\\ D\'epartement de Math\'ematiques\\BP 6759\\ F 45067 ORLEANS C\'edex 2
\\FRANCE}
\email{bonami@labomath.univ-orleans.fr\\
grellier@labomath.univ-orleans.fr \\
kacim@labomath.univ-orleans.fr}
\thanks{Authors partially supported by the 2002-2006 IHP Network, Contract Number: HPRN-CT-2002-00273 - HARP\\
The authors would like to thank Joaquim Bruna who suggested this
problem.}
\begin{abstract} We extend to  multilinear Hankel operators the fact that some truncations of  bounded Hankel operators are bounded.
We prove and use a continuity property of bilinear Hilbert
transforms on products of Lipschitz spaces and Hardy spaces.
\end{abstract}\keywords{Hankel operator, truncation, Hardy spaces, Lipschitz spaces, bilinear Hilbert transform}
\subjclass{47B35(42A50 47A63 47B10 47B49)}

\maketitle
\section{Statement of the results}

We  prove that some truncations of bounded multilinear Hankel
operators are bounded. This extends the same property for linear
Hankel operators, a result obtained by \cite{BB} (and
independently in a particular case in \cite{GB}), which we first
recall. A matrix $B=(b_{mn})_{m,n\in\N}$ is called of {\sl Hankel
type} if $b_{mn}=b_{m+n}$ for some sequence $b\in l^2(\N)$. We can
identify $B$ with an operator acting on $l^2(\N)$. Moreover, if we
identify $l^2(\N)$ with the complex Hardy space $\H2$ of the unit
disc, then $B$ can be realized as the integral operator, called
{\sl Hankel operator with symbol $b$} and denoted by $H_b$, which
acts on $f\in \H2$ by
$$H_bf(z)=
\frac{1}{2\pi}\int_{\partial\D}\frac{b(\zeta)f(\overline{\zeta})}{1-\overline{\zeta}z}d\sigma(\zeta),\;
z\in\D.$$ Here $d\sigma$ is the Euclidean measure on the unit
circle. The {\sl symbol} $b$ is given by
$b(\zeta):=\sum_{k=0}^\infty b_k\zeta^k$.  In other words,
$H_bf=\mathcal C(b\check f)$ where $\mathcal C$ denotes the Cauchy
integral, and $\check f(\zeta):=f(\overline{\zeta})$. If
$f(\zeta)=\sum_{n\in\N}a_n\zeta^n$, one has
$$H_bf(z)=\sum_{m\in\N}(\sum_{n\in\N} a_nb_{m+n})z^m.$$ Now, we
consider truncations of matrices, which are defined  as follows.
For $\beta,\gamma\in\R$, the truncated matrix
$\Pi_{\beta,\gamma}(B)$ is the matrix whose $(m,n)$ entry  is
$b_{mn}$ or zero, depending on the fact that  $m\ge \beta
n+\gamma$ or not. It is proved in \cite{BB} that such truncations,
for $\beta\neq -1$, preserve the boundedness for Hankel operators.
The proof consists in  showing that truncations are closely
related to  bilinear periodic Hilbert transforms. One then uses
the theorem of Lacey-Thiele (see \cite{LT1}, \cite{LT2},
\cite{LT3}) in the periodic setting.

We are interested in the same problem, but for multilinear Hankel
operators. For $n\in\N$, we define the multilinear Hankel operator
$H^{(n)}_b$ as follows. For $f_1,\dots,f_n \in\H2$, let
\begin{eqnarray*}H_b^{(n)}(f_1,\dots,f_n)(z)&=&\frac{1}{2\pi}
\int_{\partial\D}\frac{b(\zeta)f_1(\overline{\zeta})\dots
f_n(\overline{\zeta})}
{1-\overline{\zeta}z}d\sigma(\zeta)\cr&=&H_b(f_1\times\dots\times
f_n)(z)
\cr&=&\sum_{i_0,\dots,i_{n}}b_{m+i_1+\dots+i_n}a^1_{i_1}\dots
a^n_{i_{n}}z^m,
\end{eqnarray*}
whenever $f_j(z)=\sum a^j_{i}z^i$.  This last expression means
that, when equipped with the canonical basis of $\mathcal
H^2(\D)$, the matrix of the operator,  whose entries are indexed
in $\N^{n+1}$, are given by $b_{i_0+\dots+i_{n}}$,
${i_0,\dots,i_{n}\in\N}$. The truncated operator is obtained when
truncating this matrix as follows. For
${\mathbf\beta}=(\beta_1,\dots,\beta_n)\in \R^n$ and
$\gamma\in\R$, we denote by $\Pi_{{\mathbf
\beta},\gamma}(H_b^{(n)})$ the operator with $(n+1)$-dimensional
matrix whose $(i_0,\dots,i_{n})$-entry is $b_{i_0+\dots+i_{n}}$ if
$\beta_1 i_1+\dots+\beta_n i_n+\gamma\le i_{0}$ and zero
otherwise. Our main result is the following.

\begin{theo}\label{main}
If $H_b^{(n)}$ is a continuous multilinear Hankel operator from
$\left(\H2\right)^n$ into $\H2$, then so are its truncated
operators $\Pi_{\beta,\gamma}(H_b^{(n)})$ for any
$\beta=\nu(1,\ldots,1)$, $\gamma\in\R$. Furthermore, there exists
a constant $C=C(\nu)$ uniformly bounded when $\nu$ lies in a
compact set of $\R\setminus\{-1,0\}$ so that
$$\Vert \Pi_{\beta,\gamma}(H_b^{(n)})\Vert_{(\mathcal H^2(\D))^n\to \mathcal H^2(\D)}
\le C\Vert H_b^{(n)}\Vert_{(\mathcal H^2(\D))^n\to
\mathcal H^2(\D)}.$$
\end{theo}

Theorem \ref{main} is  deduced from an estimate on the bilinear
Hilbert transform in the periodic setting, which is of independent
interest. Let us first give some notations. The usual Lipschitz
spaces of order $\alpha$ of $2\pi$-periodic functions are denoted
by $\Lambda_\alpha(\T)$, while  $H^p(\T)$ denotes the real Hardy
space, $p>0$. Here $\T:=\R/2\pi\Z$ denotes the torus.

Let $b\in\Lambda_\alpha(\D)$ and $f\in\mathcal H^q(\D)$. The
boundary values of such functions may be seen either as functions
on the set of complex numbers of modulus 1, or on the torus
$\T=\R/2\pi\Z$, that is as $2\pi$-periodic functions on $\R$. We
write $f(e^{it})$ or $f(t)$ depending on the context, and hope
that it does not introduce any confusion.

Let us recall that, for $f$ and $b$ trigonometric polynomials on
the torus, the periodic bilinear Hilbert transform of $f$ and $b$
is given, for $k,l\in\Z$, $k\neq -l$, $l\neq 0$, by
$${\mathcal H}_{k,l} (b,f)(x)= p.v.\int_{\T} b(kx+lt)f(t) \frac{dt }{ \tan\frac{x-t} 2}.$$
Lacey-Thiele's Theorem, once transferred to the periodic setting,
is  the following.

\begin{theo}{\bf \cite{BB}}\label{BB}
Let $1<p,q\le \infty$ with $\displaystyle\frac 1r=\frac 1p+\frac 1q<\frac 32$. Then, for any $k,l\in\Z$, $k\neq -l$ and $l\neq 0$,
there exists a constant $C=C(k,l)>0$ so that, for any $f\in L^p(\T)$ and any $b\in L^q(\T)$,
$$\Vert {\mathcal H}_{k,l} (b,f)\Vert _{L^r(\T)}\le C(k,l) \Vert f\Vert _{L^p(\T)}\Vert b\Vert_{L^q(\T)}.$$
Furthermore, $C(k,l)$ is uniformly bounded when $k/l$ lies in a
compact set of $\R\setminus\{-1,0\}$.
\end{theo}

We adapt the definition to our setting, and define, for $k,l,\mu
\in \Z$ with $k\neq -l$, $l\neq 0$ and $\mu\in[-l,l]$,
\begin{equation}\label{bilinear2}
  {\mathcal  H}_{k,l,\mu}(b,f)(x)=
  \int_{\T} \left (b(kx+lt)e^{i\mu(x-t)} -b((k+l)x)\right)f((k+l)t)
\frac{dt}{ \tan{\frac{x-t} 2}}.
\end{equation}

We prove the following.
\begin{theo}\label{bilinear-p}
Let $1<p<\infty$, $0<q<p$ and  $\alpha=\frac 1q-\frac 1p$. Let $k,
l,\mu \in\Z$, with $k\neq -l$, $l\neq 0$ and $\mu\in[-l,l]$. There
exists a constant $C=C(k,l)>0$ so that, for any sufficiently
smooth functions $b\in\Lambda_\alpha(\T)$ and $f\in H^p(\T)$
\begin{equation}\label{bilin}
\Vert {\mathcal  H}_{k,l,\mu} (b,f)\Vert _{H^p(\T)}\le C\Vert f\Vert
_{H^q(\T)}\Vert b\Vert_{\Lambda_\alpha(\T)}.
\end{equation}
Furthermore, $C(k,l)$ is  independent of $\mu$ and uniformly
bounded when $k/l$ lies
 in a compact set of $\R\setminus\{-1,0\}$.
\end{theo}
In the following, we say that a constant $C(k,l)$ is {\sl
admissible} when it is uniformly bounded when $k/l$ lies in a
compact set of $\R\setminus\{-1,0\}$.

The difficulty, here, is the uniform bound of constants. We remark
that the limiting case $b\in L^\infty (\T)$ is given by the
Lacey-Thiele Theorem, that is, Theorem \ref{BB}. We will also give
a non periodic version of Theorem \ref{bilinear-p}. Our methods
rely on the ordinary Calder\`on-Zygmund theory in a bilinear
setting for the local part (one may consult \cite{GrK} and
\cite{GT}), and on transference for the non local part.

\medskip

Let us come back to holomorphic functions and to truncations.
Denote by $\Lambda_\alpha(\D)$, $\alpha>0$, the space of functions
that are holomorphic in $\D$ and whose boundary values are in
$\Lambda_\alpha(\T)$. Denote also by $\mathcal H^p(\D)$ the
complex Hardy space on the unit disc, $p>0$. Recall that, for
$p<1$, the dual of $\mathcal H^p(\D)$ is $\Lambda_\alpha(\D)$,
with $\displaystyle p=(\alpha+1)^{-1}$ (\cite{D}). As an easy
consequence of duality and factorization, one obtains that the
Hankel operator $H_b$ is bounded from $\mathcal H^q(\D)$ into
$\mathcal H^p(\D)$, with $q<p$ and $p>1$, if and only if the
symbol $b$ is in $\Lambda_\alpha(\D)$ with $\alpha=\frac 1q-\frac
1p$. Moreover, there exists a constant $C$ such that, for all
holomorphic polynomials $f$,
\begin{equation}\label{hankel}
\|H_b(f)\|_{\mathcal
H^{p}(\D)}\leq C\|b\|_{\Lambda_\alpha(\D)}\times \|f\|_{\mathcal
H^{q}(\D)}.
\end{equation}

Theorem \ref{bilinear-p} has the following corollary, which gives
the link with truncations.

\begin{coro}\label{corotroncgene}
Let $1<p<\infty$ and $0<q<p$. Let $\alpha=\frac 1q-\frac 1p$ and $b\in\Lambda_\alpha(\D)$. Then,
for any $\beta,\gamma\in\R$ with $\beta\neq -1$, the operator
 $\Pi_{\beta,\gamma}(H_b)$ is  continuous from $\mathcal H^q(\D)$
into $\mathcal H^p(\D)$. More precisely, there exists a constant $C$
uniformly bounded for $\gamma\in\R$ and $\beta$ in a compact set of
$\R\setminus\{-1,0\}$ so that
\begin{equation}\label{troncgene}
\Vert \Pi_{\beta,\gamma}(H_b)\Vert_{\mathcal H^{q}(\D)\to \mathcal
H^{p}(\D)}\le C\Vert H_b\Vert_{\mathcal H^{q}(\D)\to \mathcal
H^{p}(\D)} .
\end{equation}
\end{coro}

So, if $H_b$ is a bounded operator from $\mathcal
H^{q}(\D)$ into $\mathcal H^{p}(\D)$, its truncates
$\Pi_{\beta,\gamma}(H_b)$ are also bounded.
\begin{remark} For $\beta=-1$, we are interested in the behavior
of the  norm of the operator $\Pi_{-1,N}(H_b)$, for $N$ a positive
integer. Then $(I-\Pi_{-1,N})(H_b)=H_{S_{N-1}(b)}$ where $S_N$
gives the $N$-th-partial sum of the Fourier series, that is the
convolution operator corresponding to the Dirichlet kernel. Then,
for $1<p<\infty$, $0<q<p$, $\alpha=\frac 1q-\frac 1p$ and
$b\in\Lambda_\alpha(\D)$, one has
$$\Vert \Pi_{-1,\gamma}(H_b)\Vert_{\mathcal H^{q}(\D)\to \mathcal
H^{p}(\D)}\leq C \log N\times \Vert H_b\Vert_{\mathcal
H^{q}(\D)\to \mathcal H^{p}(\D)}.$$ Moreover, the best constant in
the previous inequality is bounded below by $c\log N$, with $c$
independent of $N$.

This result follows from the fact that $\Vert S_N\Vert_{\Lambda_\alpha\mapsto
\Lambda_\alpha}\simeq \log N$. We will prove it in the next
section for the convenience of the reader.

We use the notation $A\simeq B$ whenever there exist two universal
constants $c,c'>0$ so that $c'B\le A\le cB$.

\end{remark}
\begin{remark} \label{remarkbeta}For $\beta=0$, it is easy to show that truncations preserve uniformly
the class of bounded Hankel operators. It is sufficient to
consider positive integer values of $\gamma$. In fact, one has
$\Pi_{0,N+1}(H_bf)=(I-S_{N})H_bf$ for $N\in \N$ so that it
suffices to use that $S_N$ is a bounded operator on $\mathcal
H^p(\D)$, $1<p<\infty$, with a bound independent on $N$. When
$\beta=\infty$ (that is when truncations of Hankel matrices are
with respect to directions $n=N$),
 one may show also the analogous fact. To do so, it suffices to argue by duality and to
 show that, for $b\in\Lambda_\alpha(\D)$, $\Pi_{0,N}(H_b)$ maps boundedly $\mathcal H^{p'}(\D)$
 into $\Lambda_{\alpha'}(\D)$, where $\alpha'=\alpha-1/p'$, with a bound independent of $N$.
 We will prove this in the next section.
\end{remark}

Let us deduce Theorem \ref{main} from the corollary. It is clear
that $$H_b^{(n)}(f_1,\dots,f_n)(z)=H_b(f_1\times\dots\times
f_n)(z).$$ Using the factorization of functions in Hardy classes,
we know that $H_b^{(n)}$ is bounded as an operator from
$\left(\H2\right)^n$ into $\H2$ if and only if the Hankel operator
$H_b$ is bounded  from $\mathcal H^{2/n}(\D)$ into $\H2$, that is,
if and only if $b$ is in $\Lambda_\alpha(\D)$ for $\alpha=\frac
{n-1}2$ . To conclude, we use the fact that the truncation
$\Pi_{{\beta},\gamma}$ of $H_b^{(n)}$ for ${
\beta}=\nu(1,\dots,1)$, corresponds to the truncation
$\Pi_{\nu,\gamma}$ of $H_b$, as it can be easily verified.

\bigskip
Let us finally remark that this paper leaves open the general
problem of all truncations of multilinear Hankel operators. One
may be tempted to reduce to multilinear Hilbert transforms, as
studied in \cite{MTT}.
\bigskip

 The remainder of the paper is organized as
follows. In the next section, we deduce the corollary from Theorem
\ref{bilinear-p}. In the last one, we prove Theorem \ref{main}.
Let us emphasize the fact that this last proof does not use
Lacey-Thiele Theorem, and is elementary compared to it.

\section {The link between truncations and bilinear Hilbert
transforms} At first, we prove the result corresponding to
truncations with $\beta=\infty$, as stated in the remark
\ref{remarkbeta}.
\begin{proof}
Let us recall that it is sufficient to prove that, for
$\alpha'=\alpha-\frac 1{p'}$,
$$\Vert \Pi_{0,N}H_b(f)\Vert_{\Lambda_\alpha'(\D)}\le C \Vert f
\Vert_{\mathcal H^{p'}(\D)}\Vert b\Vert_{\Lambda_\alpha(\D)}$$
with a constant independent of $N\in\N$. By density, it is
sufficient to prove this inequality on trigonometric polynomials.
It is easy to show that
$$Tf(z):=\Pi_{0,N}H_b(f)(z)= \frac 1{2\pi}\int_{\partial\D} b(\zeta)\check
f(\zeta)(z\overline\zeta)^{N}\frac{d\sigma(\zeta)}{1-z\overline{\zeta}},$$
where $\check f(\zeta):=f(\overline \zeta)$ as before. Using
Stokes Formula, we get
$$Tf(z)=\frac 1{\pi}\int_{\D}
\left(b(re^{i\theta})+b'(re^{i\theta})\right)f(re^{-i\theta})z^Nr^Ne^{-iN\theta}\frac{d\theta
dr}{1-zre^{-i\theta}}.$$

Assume to simplify that $0<\alpha<1$. By the holomorphic
characterization of Lipschitz spaces, it suffices to show that
$$\vert (Tf)'(z)\vert \le C(1-\vert z\vert)^{\alpha'-1}.$$
Computing the derivative of $T(f)$, and using the fact that $\vert
Nz^Nr^N\vert\le  C{(1-r|z|)^{-1}}$, we get the estimate
$$
|(Tf)'(z)|\leq C \int_{\D} \left |
b(re^{i\theta})+b'(re^{i\theta})\right|\,|f(re^{-i\theta})|(1-r|z|)^{-1}\frac{d\theta
dr}{|1-zre^{-i\theta}|}.
$$
We use that $\vert b(re^{i\theta})+b'(re^{i\theta})\vert \le C
\Vert b\Vert_{\Lambda_\alpha(\D)}(1-r)^{\alpha-1}$  and H\"older's
Inequality to get that
\begin{eqnarray*}
|(Tf)'(z)|&\le&C\Vert b\Vert_{\Lambda_\alpha(\D)}\Vert
f\Vert_{L^{p'}(\T)}\int_{0}^1
\frac{(1-r)^{\alpha-1}}{(1-r|z|)^{2-1/p}}dr\cr &\le&
C(1-|z|)^{\alpha-1+1/p}\Vert b\Vert_{\Lambda_\alpha(\D)}\Vert f
\Vert_{\mathcal H^{p'}(\D)}.
\end{eqnarray*}
It ends the proof.
\end{proof}

We now recall the following characterization of Lipschitz
spaces.
\begin{proposition}\label{propLip}
A function $b$ belongs to $\Lambda_\alpha(\T)$ if and only if
there exists a decomposition $b=\sum_{j=0}^\infty b_j$ with
$\Vert\nabla^l b_j\Vert_{\infty}\le A2^{-j(\alpha-l)}$ for any
integer $0\le l\le k$, where $k$ is the smallest integer
$>\alpha$.

If $b\in\Lambda_\alpha(\D)$, the decomposition may be chosen so
that the spectrum of $b_j$, $j\ge 1$, lies in
$\{2^{j-1},\cdots,2^{j+2}\}$ and $b_0$ is a constant. Moreover,
the norm is equivalent to the smallest constant $A$ for which we
have $\Vert b_j\Vert_{\infty}\le A2^{-j\alpha}$ for all $j\ge 0$.
\end{proposition}

The first statement can be found in \cite{S}, on page 256. The
second one can be obtained from the proof given there when
$\alpha<1$, using a variant of the kernel of De La Vall\'{e}e Poussin
adapted to the dyadic decomposition. Then one concludes using the
action of powers of the Laplacian, which give  isomorphisms
between Lipschitz spaces with different parameters.
\medskip

As a classical corollary, we prove the following.
\begin{coro} For any $N\in\N$, the norm of the operator $S_N$ from
$\Lambda_\alpha(\T)$ into itself is equivalent to $\log N$. In other
words $$\Vert S_N\Vert_{\Lambda_\alpha(\T)\mapsto \Lambda_\alpha(\T)}\simeq \log
N.$$
\end{coro}
\begin{proof}
Denote by $D_N$ the Dirichlet kernel of order $N$. As $\Vert D_N\Vert_{L^1(\T)}\simeq \log N$,
it follows easily that $$\Vert S_N \Vert_{\Lambda_\alpha(\T)\mapsto \Lambda_\alpha(\T)}\leq C\log N.$$
It is sufficient to prove the converse inequality for the operator
$S_N^+$ defined by $$S_N^+f(t)=\sum_{n\leq N} \widehat f (n)
e^{int}.$$ We start from the well known fact that, for all $M\geq 1$, one can find a
trigonometric polynomial $g$ of degree $2M$ such that
$\|g\|_{L^\infty(\T)} \simeq 1$ and $\|S_M^+g\|_{L^\infty(\T)} \simeq \log M$.  Let
us assume that $N=4\times 2^n$ for some integer $n\geq 1$. Take $M=\frac
N4$, then the function $f=e^{i3\frac N4t}g$ satisfies $\|f\|_{L^\infty(\T)} \simeq 1$
and $\|S_N^+f\|_{L^\infty(\T)} \simeq \log N$ and has a spectrum  included in $\{2^{n},\cdots, 2^{n+3}\}$. Moreover, by the previous proposition, and the spectrum
assumption on $f$, one has $\Vert f\Vert_{\Lambda_\alpha(\T)}\simeq
2^{n\alpha}\Vert f\Vert_{L^\infty(\T)}$ and $\Vert
S_N^+f\Vert_{\Lambda_\alpha(\T)}\simeq 2^{n\alpha}\Vert S_N^+
f\Vert_{L^\infty(\T)}$. This finishes the proof.
\end{proof}

In the following, we are going to use the same kind of idea in the
general case to get rid of the dependence in $\gamma$ of the
constants. Let us consider now the truncation
$\Pi_{\beta,\gamma}$. Remark that, if $\beta$ and $\gamma$ are
negative, then $\Pi_{ \beta,\gamma}H_b=H_b$ and there is nothing
to prove. We will not consider this case any longer. If $\gamma<0$
and $\beta>0$, then
$(I-\Pi_{\beta,\gamma})H_b=(I-\Pi_{\beta,\gamma})H_{P_{[-\gamma/\beta]}(b)}$
where $P_N(b)$ denotes the projection of $b$ on the space of
functions whose spectrum is contained in $[N,\infty)$. If
$\gamma>0$ and $\beta>0$, then
$\Pi_{\beta,\gamma}H_b=\Pi_{\beta,\gamma}H_{P_{[\gamma]}(b)}$. In
the last case, where $\gamma>0$ and $\beta<0$, then
$\Pi_{\beta,\gamma}H_b=\Pi_{\beta,\gamma}H_{P_{n^{\beta,\gamma}}(b)}$
where $n^{\beta,\gamma}=\min\{[\gamma],[-\gamma/\beta]\}$ (here
$[.]$ denotes the integer part). So, in any case, $b$ may be
replaced by some projection $P_{n^{\beta,\gamma}}(b)$ of $b$ where
$$n^{\beta,\gamma}=\begin{cases}
=[-\gamma/\beta]\text{ when }\gamma<0\text{ and }\beta>0\cr
=[\gamma]\text{ when }\gamma>0\text{ and }\beta>0\cr
=\min\{[\gamma],[-\gamma/\beta]\}\}\text{ otherwise
}\end{cases}.$$ So, it suffices to consider the symbols $b$ modulo
a polynomial of degree smaller than  $n^{\beta,\gamma}$. The
following lemma will allow us to modify $b$ by adding such a
polynomial. The main point, here, is the fact that the constant is
 admissible.\begin{lem}
\label{lemmeLipschitz} Let $\alpha>0$. There exists a constant $C$
with the following property: for any $b\in\Lambda_\alpha(\D)$ and
$N\in\N$, there exists a function $\tilde {b}$ such that $P_N
(\tilde b)=P_{N}(b)$, and, for all $M\in\Z$,
\begin{eqnarray*}
\Vert \tilde b{\zeta}^{M}\Vert_{\Lambda_\alpha(\T)}&\le& C\left(\frac {|M|}
{N+1}+1\right)^\alpha\Vert b\Vert_{\Lambda_\alpha(\D)}.
\end{eqnarray*}
\end{lem}
\begin{proof} When $N\leq 16$, we take $\tilde b =b$. We conclude easily, using the fact
 that the  norm in $\Lambda_\alpha(\T)$ of $\zeta\to
{\zeta}^{M}$ is equivalent to $(|M|+1)^\alpha$.  Consider now
$N>16$, and choose $N_0$ so that $2^{N_0}\le N <2^{N_0+1}$. Use
the Littlewood-Paley decomposition given in Proposition
\ref{propLip} to write $b$ as $b=\sum b_j$ with the spectrum of
$b_j$, for $j\ge 1$, included in $[2^{j-1},2^{j+2})$. We can
clearly take $$\tilde b = \sum_{j\ge N_0-2}b_j\ .$$  For
simplicity, we assume that $\alpha<1$. The new function $a(\zeta)=
\tilde b (\zeta){\zeta}^{M}$ can be written as $\sum_{j\ge
N_0-2}a_j$, with $a_j(\zeta)=b_j(\zeta){\zeta}^{M}$. It is clear
that, for $j\ge N_0-2$,
\begin{equation}\label{bernstein}
  ||a_j||_{\infty}\le C 2^{-j\alpha}, \ \ \ \ \ \ ||a_j'||_{\infty}\le C \left(\frac
{|M|} {N+1}+1\right) 2^{j(1-\alpha)}\ .
\end{equation}
 Indeed, for the
second inequality, we use Bernstein Inequality to get that
$||b_j'||_{\infty}\leq C 2^j||b_j||_{\infty}$. The second term, in
the derivative of $a_j$, is bounded by
$|M|\times||b_j||_{\infty}$. The inequality follows at once, using
the fact that $N+1\leq  2^{j+3}$. We adapt the proof given in
\cite{S} to deduce the required inequality for $a$.

\end{proof}
\begin{proof}

We are now in position to prove that Corollary \ref{corotroncgene} follows from Theorem \ref{bilinear-p}.
It is sufficient to prove, for trigonometrical polynomials $b$ and $f$, that
\begin{equation}\label{ineq}
\|\Pi_{ \beta,\gamma}H_b(f)\|_{\mathcal H^p(\D)}\leq
C\|b\|_{\Lambda_\alpha(\D)}\times \|f\|_{\mathcal H^{q}(\D)}
\end{equation}
for some constant $C$ that depends only on $\beta$ and is
uniformly bounded when
 $\beta$ lies in a compact set of $\R\setminus \{-1,0\}$.
Moreover, since $b$ and $f$ are trigonometric polynomials, we can restrict to rational values $\frac kl$ for
$\beta$, and values of $\gamma$ such that $\gamma l$ is an
integer, as soon as we show that the constant $C$ may be taken
independent of $k,l, \gamma$ when these last ones vary in such a
way that $\frac kl$ lies in a compact set of $\R\setminus
\{-1,0\}$. So, we assume that $\displaystyle\beta=\frac
kl\neq -1$ with $k,l\in\Z$, $l\neq 0$ and that $l\gamma\in\Z$. We
may also assume that $k+l>0$.

Assume that $f(t)=\sum_{n\in\N} a_n e^{int}$. An elementary
computation shows that the analytic part of the function
$$ p.v.\int _\T b(kx+lt)e^{i\gamma l(x-t)} \check{f}((k+l)t)
\frac{dt}{\tan \frac{x-t}{2}}$$ is equal to
$$\sum_{n,m\in\N} a_n b_{m+n} \text{sign}(l(m-\beta
n-\gamma))e^{im(k+l)x}. $$ We recognize the function $
\text{sign}(l)(2\Pi_{\beta,\gamma}-I)(H_b(f))(e^{i(k+l).})$, for
which we want to have an $L^p$ estimate (recall that we have such
an estimate for  $H_b(f)$). By Lemma \ref{lemmeLipschitz} and the
remark before, we can replace $b$ by the function $\tilde b$ that
corresponds to the choice $N=n^{\beta, \gamma}$. On the other
hand, this term is equal to the analytic part of the function
$$ e^{i[\gamma](l+k)x}\int _\T \underbrace{\tilde b(kx+lt)e^{-i[\gamma] (kx+lt)}}_{=\hat b(kx+lt)} e^{i\delta l(x-t)}\check{f}((k+l)t
\frac{dt}{\tan \frac{x-t}{2}}$$ where $\delta=\gamma-[\gamma]$.
This is the analytic part of $e^{i[\gamma](l+k)x}{\mathcal
H}_{k,l,\delta l} (\hat b,f)$, up to the Cauchy projection of
$$\tilde b((k+l)x)\times p.v.\int _\T \check{f}((k+l)t)
\frac{dt}{\tan \frac{x-t}{2}}.$$ This last term coincides with
$H_{\tilde b}(f)(e^{i(k+l)x})$ and is in $\mathcal H^p(\D)$ with
norm bounded by
$$C\|\tilde b\|_{\Lambda_\alpha(\D)}\times
\|f\|_{\mathcal H^{q}(\D)} \le C\|b\|_{\Lambda_\alpha(\D)}\times
\|f\|_{\mathcal H^{q}(\D)}. $$ By Theorem \ref{bilinear-p}, the
term $e^{i[\gamma](l+k)x}{\mathcal  H}_{k,l,\delta} (\hat
b,f)$ is also in $H^p(\T)$ with norm bounded by $C\Vert \hat
b\Vert_{\Lambda_\alpha(\T)}\Vert f\Vert_{ H^{q}(\T)}$  for some
admissible constant $C$. This last quantity is, in turn, bounded
by  $C\left(|\beta|+1\right)^\alpha
\|b\|_{\Lambda_\alpha(\D)}\times \|f\|_{\mathcal H^{q}(\D)}$, as
one can see by using Lemma \ref{lemmeLipschitz}.
\end{proof}
\section{Proof of the Theorem \ref{bilinear-p}.}

Let us note $L=k+l$. Let us remark first that it is sufficient to
consider the case when $L>0$ (the general case  follows by
replacing $b$ by $\tilde b$ such that $\tilde b(x):=b(-x)$ and $f$
by $\tilde f$ such that $\tilde f(x):=f(-x)$).

Let us give a first  reduction of the problem. We use the
following lemma, where $\tau$ denotes the translation by $y$ on
the torus. Its proof is elementary.
\begin{lem} For all trigonometrical polynomials $b$ and $f$,  the
following identity is valid.
\begin{equation}\label{transl}
{\mathcal H}_{k,l,\mu}(b\circ \tau^L,f \circ \tau )= {\mathcal
H}_{k,l,\mu}(b,f)\circ\tau.
\end{equation}
\end{lem}

 Cutting into eight parts the function $f$ and using this lemma, we reduce
 to periodic functions with support in $(0, +\frac \pi4)$
in the case $q>1$. For the case of $q\leq 1$, by the atomic
decomposition Theorem of $H^q(\T)$, it suffices to consider the
action of ${\mathcal H}_{k,l,\mu}(b,\cdot)$ on $H^q(\T)$-atoms,
also called $q$-atoms. So, let $f$ be a $q$-atom, that is,  either
a constant or a function $f$ on $\T$ supported in some interval
$I$ on the torus (we will always assume its length less than
$\pi/4$, which is possible)  so that
$$\Vert f\Vert_{L^\infty(\T)}\le |I|^{-1/q} \text{ and }\int_\T x^kf(x)dx=0\text{ for any integer }k\le \frac{1}q-1.$$
% If $f$ is the constant atom, or if $f$
%is a non constant atom but whose support is  some interval $I$ of
%the torus of length bigger than $\frac \pi4$, then, for all $r>1$,
%its $L^r$ norm is uniformly bounded. It follows at once that the
%$L^p$ norm of ${\mathcal H}_{k,l,\mu}(b,f)$ is also uniformly
%bounded, assuming that we have already proved the required
%inequality in $H^r(\T)$, for $r>1$. So we can restrict to atoms $a$ whose
%support $I$ has length less than $\frac \pi4$. Using (\ref{transl}), we can assume that $I=(0,r)$, as soon as constants
%depend only on the Lipschitz norm of $b$, which is invariant by
%translation.

  We have seen that
${\mathcal H}_{k,l,\mu}(b,f)$ is $2\pi/L$-periodic. We simplify
the notation, and write
$$Tf(x):={\mathcal H}_{k,l,\mu}(b,f)\left(\frac xL\right),$$
so that, after a change of variable in the integral, one can write
$$Tf(x)=\int_{-\pi L}^{\pi L} \left[ b\left(x+\frac{l}{L}(t-x)\right)e^{i\frac{\mu}{L}(x-t)}
-b(x)\right]f(t)\frac{dt}{\tan{\frac{x-t} {2L}}}.$$ We are
restricted to prove that
\begin{equation}\label{final}
\Vert T(f)\Vert_{L^p(]-\pi,\pi[)}\leq C \Vert
b\Vert_{\Lambda_\alpha(\T)}\Vert f\Vert_{H^q(\T)}
\end{equation}
for $f$ supported in the interval $(0, +\frac \pi4)$ when $q>1$,
or
\begin{equation}\label{final-bis}
\Vert T(f)\Vert_{L^p(]-\pi,\pi[)}\leq C \Vert
b\Vert_{\Lambda_\alpha(\T)}\Vert f\Vert_{H^q(\T)}
\end{equation}
for $f$ is a $q$-atom supported in $(0, +r)$ for some $r\leq \frac
\pi4$ when $q\leq 1$.

\medskip
 Let us explain the next reduction of the problem.
 For $|x|<\pi$,
 $|s|<L\pi$ and $s$ in the support of the periodic function $f$, we can
write
$$\frac 1{L\tan\frac{x-s}{ 2L}}=\frac 1{x-s}
+ \frac 1L\psi \left(\frac{x-s}L\right),$$ with $\psi$ a ${\mathcal
C}^\infty$ function with compact support. Indeed, the conditions on $s$, $x$ and $L$, imply that
$\frac{|x-s|}L\leq \frac {3\pi}2$. We get rid of
the term with $\psi$ by showing the next lemma.

\begin{lem}
Let $\psi$ be a ${\mathcal C}^\infty$ function with compact
support. Then there exists admissible constants $C, C_q$ depending
on $\psi$ such that, for $f$ a periodic integrable function supported in $(0, +\pi/4)$ (modulo $2\pi$), the quantity
$$A(x)=\frac 1L\int_{-\pi L}^{\pi L} \left[ b\left(x+\frac{l}{L}(s-x)\right)
e^{i\frac{\mu}{L}(x-s)}-b(x)\right]f(s)\psi
\left(\frac{x-s}L\right)ds $$ is uniformly bounded, for $|x|<\pi$,
by $C \Vert b\Vert_{L^\infty(\T)}\|f\|_{L^1(\T)}$. Moreover, it is bounded by $C_q\Vert b\Vert_{\Lambda_\beta(\T)}$ when $f$ is a
$q$-atom, with $\beta=1/q-1$.
\end{lem}
\begin{proof}
  $A(x)$ may be
written as a scalar product $\frac 1L <\tilde f,g>$ on the real
line, where $\tilde f$ is the function on the real line which
coincides with $f$ on $(-L\pi, +L\pi)$ and is zero outside, while
 $g$ takes care of the other terms. Consider the first case, and remark
 that $\Vert \tilde f\Vert_{L^1(\R)}= L\Vert f\Vert_{L^1(\T)}$, while the norm
of $g$ in $L^\infty(\R)$ is bounded, up to a constant, by the norm
of $b$ in $L^\infty(\T)$. The conclusion follows at once. Assume
now that $q<1$. Then $\tilde f$ is the sum of $L$ $q$-atoms of
$H^q(\R)$. On the other hand, one has $\Vert
g\Vert_{\Lambda_\beta(\R)}\le C\Vert b\Vert_{\Lambda_\beta(\T)}$,
with $C$ admissible constant. We conclude, using the fact that the
space $\Lambda_\beta(\R)$ being the dual space of $ H^q(\R)$.
\end{proof}
\medskip

The main term can be  written as $\widetilde
Tf(x)+\overline{T}f(x)$, where
\begin{equation}\label{tilde}
\widetilde Tf(x)=\int_{-\pi}^{\pi} \left[
b\left(x+\frac{l}{L}(s-x)\right)
e^{i\frac{\mu}{L}(x-s)}-b(x)\right]f(s)\frac{ds}{x-s}
\end{equation}
and
\begin{equation}\label{bar}
\overline Tf(x)=\int_{\pi<|s|<L\pi} \left[
b\left(x+\frac{l}{L}(s-x)\right)
e^{i\frac{\mu}{L}(x-s)}-b(x)\right]f(s)\frac{ds}{x-s}.
\end{equation}
Let us start with the first term.

\begin{prop}\label{prop1}
There exist admissible constants
$C_q$ so that
\begin{equation}\label{homog}
\Vert\widetilde T(f)\Vert _{L^p(-\pi, +\pi)}\le C_q \Vert
b\Vert_{\Lambda_\alpha(\T)}
\end{equation}
when $q>1$ and  $f$ is supported in the interval $(0, +\frac
\pi4)$ and $\Vert f\Vert _{L^q(\T)}\leq 1$, or when $q\leq 1$ and $f$ is a
q-atom supported in $(0, +r)$, for some $r\leq \frac \pi4$.
\end{prop}

\begin{proof}
We can deduce the required estimate from the analogue of Theorem
\ref{bilinear-p} on the real line, which we state now. We define,
for  $\beta\in \R$, for $b\in\Lambda_\alpha(\R)$ and $f\in
H^q(\R)$,
\begin{equation}\label{real-op}
  {\mathcal  H}_{\beta} (b,f)(x):=\int_{-\infty}^{\infty} \left[ b\left(x+\beta(s-x)\right)
-b(x)\right]f(s)\frac{ds}{x-s} .
\end{equation}
We have the following proposition.
\begin{prop}\label{real}
Let $1<p<\infty$, $0<q<p$ and  $\alpha=\frac 1q-\frac 1p$. There exists a constant $C=C(\beta)>0$ so that, for
any sufficiently smooth functions $b\in\Lambda_\alpha(\R)$ and
$f\in H^q(\R)$
\begin{equation}\label{bilin2}
\Vert {\mathcal  H}_{\beta} (b,f)\Vert _{H^p(\R)}\le C\Vert f\Vert
_{H^q(\R)}\Vert b\Vert_{\Lambda_\alpha(\R)}.
\end{equation}
 Furthermore, $C(\beta )$ is uniformly bounded when $\beta$ lies
 in a compact set of $\R$.
\end{prop}
\begin{remark} Let us emphasize that this result may be seen as an elementary case of
Lacey-Thiele Theorem on the bilinear Hilbert transform (see
\cite{LT1}, \cite{LT2}, \cite{LT3}).
\end{remark}

\begin{proof}
 Let us again simplify the notation, and set $S(f):= {\mathcal  H}_{\beta} (b,f)$. The kernel of $
S$ is bounded, up to the constant $\displaystyle
C|\beta|^\alpha\Vert b\Vert_{\Lambda_\alpha(\R)}$, by the Riesz potential
$|x-s|^{-1+\alpha}$, and the estimate follows directly from the
classical estimates on Riesz potentials when $q>1$.

Let us now concentrate on $q\leq 1$, and assume that $f$ is an
atom supported in the interval $I =(-r,+r)$ (we can reduce to this
case, using invariance by translation as in Lemma \ref{transl}).
Using interpolation, it is sufficient to consider non integer
values of $\alpha$. We assume that $k<\alpha<k+1$. We write $
S(f)=A_1+A_2 $, with $A_1=S(f)1\! \text{l}_{(-2r,+2r) }$. We prove
that both $A_1$ et $A_2$ are $L^p$-functions. To prove that  $A_1
\in L^p$, we remark that, because of the fact that $f$ has
vanishing moments up to order $k$, we can replace the content of
the bracket, in the definition of $S(f)$ as an integral given in
(\ref{real-op}), by
$$ b\left(x+\beta(s-x)\right)
-\sum_{j=0}^k \frac {b^{(j)}(x)}{j!}\beta^j(s-x)^j.$$ This last
quantity is bounded by $C |\beta|^\alpha |s-x|^\alpha\Vert
b\Vert_\alpha$. We obtain that
$$\vert A_1(x)\vert\le C |\beta|^\alpha r^{\alpha-\alpha'}\Vert
b\Vert_\alpha \mathcal I_{\alpha'}(|f|)(x),$$ where $\mathcal
I_{\alpha'}$ denotes the fractional integral related to the Riesz
potential $|x-s|^{-1+\alpha'}$. We have chosen $\alpha'<1$, such
that $\frac 1p=\frac1{\tilde p}-\alpha'$ with $\tilde p>1$. Then
$$\Vert A_1\Vert_{L^p}\le  C |\beta|^\alpha r^{\alpha-\alpha'}\Vert
b\Vert_\alpha \Vert f\Vert_{L^{\tilde p}}.$$ We conclude for $A_1$
using the fact that $f$ is a  $q$-atom supported in $(-r,+r)$,
that is, $\Vert f\Vert_{L^\infty(\R)}\leq C r^{-1/q}$.

To deal with $A_2$, we replace the content of the bracket of
(\ref{real-op}) differently. We write that $
b\left(x+\beta(s-x)\right)-b(x)$ is the sum of three terms, that
is,
  \begin{eqnarray*}
      b\left((1-\beta)x+\beta s)\right)
-\sum_{j=0}^k \frac {b^{(j)}((1-\beta)x)}{j!}\beta^j s^j \hspace {3cm} \\
     +\left (\sum_{j=0}^k \frac {b^{(j)}((1-\beta)x)}{j!}\beta^j x^j
     -b(x)\right )+\sum_{j=1}^k \frac {b^{(j)}((1-\beta)x)}{j!}\beta^j
     (s^j-x^j).
  \end{eqnarray*}
  The last term may be written as the product of $s-x$ with a
  polynomial in $s$ of degree less than $k$. Because of the
  condition on the moments of $f$, the corresponding term is zero.
  The first term is bounded by $ r^\alpha |\beta|^\alpha\Vert
b\Vert_{\Lambda_\alpha(\R)}$, the second one by $  \vert\beta\vert^\alpha\Vert
b\Vert_{\Lambda_\alpha(\R)}|x|^\alpha$. So, we can write $A_2$ as a sum of two
terms, $A_2^{(1)}$ and
  $A_2^{(2)}$, corresponding to the two terms above.
For $|x|>2r$, $\vert A_2^{(1)}(x)\vert $ is bounded by $ r^\alpha |\beta|^\alpha\Vert
b\Vert_{\Lambda_\alpha(\R)} \Vert f\Vert_{L^{ 1}(\R)}/|x|$. The bound for its $L^p$
norm follows at once. For $A_2^{(2)}$, we use the well-known fact
that the Hilbert transform of the $q$-atom $f$ is bounded by
$Cr^{k+2-1/q}/|x|^{k+2}$ for $|x|>2r$. We obtain the estimate
$$|A_2^{(2)}(x)|\leq C  |\beta|^\alpha\Vert b\Vert_{\Lambda_\alpha(\R)} r^{k+2-1/q}|x|^{\alpha-k-2}$$
for $|x|>2r$, from which we conclude at once.

This finishes the proof of Proposition \ref{real}.
\end{proof}
We conclude for the proof of Proposition \ref{prop1} using
Proposition \ref{real} for ${\mathcal  H}_{\beta} (\tilde b,\tilde
f)$, with $\beta=\frac lL$, $\tilde f$ the function that is
equal to $f$ on $(0,r)$ and vanishes elsewhere, and $\tilde b(x) =
b(x)e^{-i\frac \mu l x}$ (it is elementary to show that $\tilde b$
belongs to $\Lambda_\alpha(\R)$ with a Lipschitz norm bounded by
$C\Vert b\Vert_{\Lambda_\alpha(\T}$).
\end{proof}

We now turn to the estimate of $\overline{T}f$. We prove the following.
\begin{prop}
Let $p>1$. There exist admissible constants $C,C_q$ so that,
\begin{enumerate}
\item if $f\in L^1(\T)$ and is  supported in the interval $(0,
+\frac \pi4)$, then
\begin{equation}\label{bar-1}
\Vert\overline T(f)\Vert _{L^p(-\pi, +\pi)}\le C\Vert f\Vert
_{L^1(\T)}\Vert b\Vert_{{L^p(\T)}}.
\end{equation}
\item if $q<1$  and
$f$  is a $q$-atom supported in $(0, +r)$, for some $r\leq \frac
\pi4$,  then
\begin{equation}\label{bar-at}
\Vert\overline T(f)\Vert _{L^p(-\pi, +\pi)}\le C_q\Vert
b\Vert_{\Lambda_\beta(\T)},
\end{equation}
with $\beta= 1/q-1$.
\end{enumerate}
\end{prop}
\begin{proof}
We may assume that $L\geq 2$, otherwise there is nothing to prove.
We begin with the case $q\ge 1$. We will not use any compensation
between the two terms involving $b$, so we replace $\overline T$
by $S$, where
$$S(f)(x):=\int_{\pi<|t|<\pi L} b\left(\frac{k}{L}x
+\frac{l}{L}t\right)e^{i\frac \mu L (x-t)} f(t)\frac{dt}{x-t}.$$
We  show (\ref{bar-1}) for $S(f)$ in place of $\overline T(f)$.
The difference between both is of the same type as $S(f)$, except
that $b$ is replaced by a constant. Now, we use the
$2\pi$-periodicity of $f$ to write $S(f)(x)$ as
$$\int_{-\pi}^{\pi} f(t)\sum_{1\le |j|< L/2}b\left(\frac{k}{L}x
+\frac{l}{L} (t+2\pi j)\right)e^{i\frac \mu L (x-t-2\pi
j)}\frac{dt}{x -2\pi j-t}.$$ We may replace ${t +2\pi j-x}$ by
$2\pi j$ in the last integral. Indeed, the difference is bounded,
up to a constant, by $1/j^2$. Taking the sum,  we get a quantity
which is bounded, up to an admissible constant, by $\Vert b\Vert
_{L^p(\T)}\Vert f\Vert_{L^1(\T)}$. Indeed, the constant comes from
the fact that,  for the periodic function $b$,
\begin{equation}\label{compar}
  \left\Vert b(\frac kL\cdot)\right\Vert_{L^p(-\pi, +\pi)}
\le C \left(1+\frac L{|k|}\right)^{\frac 1p} \Vert
b\Vert_{L^p(\T)}.
\end{equation}

It remains to consider the term
$$\int_{-\pi}^{\pi} f(t)\sum_{1\le |j|< L/2} b\left(\frac{k}{L}x
+\frac{l}{L} (t+2\pi j)\right) \frac{ e^{-i\frac \mu L (t+2\pi j)}}{j}dt.$$ Taking
the $L^p$-norm on $(-\pi, +\pi)$ and using Minkowski inequality,
it is sufficient to show that $\Vert B(\frac kL\cdot, t)\Vert
_{L^p(-\pi, +\pi)}\le C\Vert b\Vert_{L^p(\T)}$, with $C$ an admissible
constant which is independent of $t$, where
$$ B(x,t) : = \sum_{1\le |j|< L/2}b\left(x
+\frac{l}{L} (t+2\pi j)\right) \frac{e^{-2\pi i\frac {\mu j} L }}j.$$ We
 use the inequality (\ref{compar}) written for the function
$B(\cdot, t)$ to reduce to an inequality for this last function.
We will use transference to conclude. Let $\Theta_j$ be the
strongly continuous representation of $\Z$ in $L^p(\T)$ given by
$$\Theta_jb(t)= b\left(t+\frac{l}{L}2\pi j\right).$$ It is clear that
$B(\cdot, t)$ is obtained from  $b(\cdot+\frac{lt}{L})  $ when
using the operator  $$\sum_{1\le |j|< L/2}\frac{e^{-2\pi i\frac
{\mu j} L }}j \Theta_j.$$ The theory of tranference (see \cite{CW})
allows us to conclude that it is a uniformly bounded operator on
$L^p(\T)$ once we know that the operator given by the convolution
by the sequence $c_j:= e^{-2\pi i\frac {\mu j} L }/j$ for $1\le
|j|< L/2$, $c_j=0$ otherwise, is bounded on $\ell^p(\Z)$. This,
again, is classical: the convolution operator on $\ell^p(\Z)$
associated to the sequence $(1/j)_{|j|\ge 1}$ is a bounded
operator on $\ell^p(\Z)$, $p>1$, as the discrete analogue of the
Hilbert transform. The same holds also for the truncated sequence
$\left(c_j^L\right)$ where $c_j^L:= 1/j$ for $1\le |j|< L/2$, and
$c_j^L=0$ otherwise, with an operator norm bounded independently
of $L$. The same is also valid for the sequence $c_j$ given above:
we may write this last operator in terms of the conjugate of the
previous one under the action of the multiplication by $e^{-2\pi
i\frac {\mu j} L }$. Finally,
$$\left\Vert B(\cdot, t)\right\Vert _{L^p(\T)}\le C\left\Vert
b(\cdot+\frac{lt}{L}) \right\Vert_{L^p(\T)}= C\Vert b\Vert_{L^p(\T)}.
$$
 This finishes the proof of (\ref{bar-1}).

\medskip

Let us now prove (\ref{bar-at}). For simplicity, we  only consider
the case when $\beta$ is not an integer, which is sufficient for
our purpose.  We assume that $f$ is an atom supported in $(0,r)$.
We write $\overline T(f)(x)$ as
$$\int_{-\pi}^{\pi} f(t)\sum_{1\le |j|< L/2}\left[b\left(\frac{k}{L}x
+\frac{l}{L} (t+2\pi j)\right)e^{i\frac \mu L (x-t-2\pi
j)}-b(x)\right]\frac{dt}{x -2\pi j-t}.$$ As in the proof of
Proposition \ref{real}, we  write the content of the bracket as a
sum of terms. Let us first define
\begin{equation}\label{new-b}
  \tilde b(x):=b(x)e^{-i\frac \mu l x}.
\end{equation}
As we  already pointed out, it is elementary to see that $\tilde
b$ is in $\Lambda _\beta(\R)$ with a norm bounded by $C\Vert
b\Vert_{\Lambda_\beta(\T)}$. Then, the content of the bracket is equal
to $e^{i\frac \mu l x}\tilde B_j(x, t)$ with $\tilde B_j(x,
t)=\tilde B_j^{(1)}(x, t)+\tilde B_j^{(2)}(x, t)+\tilde
B_j^{(3)}(x, t)$, with
\begin{eqnarray*}\label{formule}
 \tilde B_j^{(1)}(x, t)&=& \tilde b\left(\frac{k}{L}x
+\frac{l}{L} (t+2\pi j)\right)-\sum_{s=0}^{[\beta]}\tilde
b^{(s)}\left(\frac{k}{L}x +\frac{l}{L}2\pi j\right)\left(\frac l L\right)^s
\frac{t^s}{s!}\cr \tilde B_j^{(2)}(x,
t)&=&\sum_{s=0}^{[\beta]}\tilde b^{(s)}\left(\frac{k}{L}x
+\frac{l}{L}2\pi j\right)\left(\frac l L\right)^s \frac{(x-2\pi
j)^s}{s!}-\tilde b(x)\cr \tilde B_j^{(3)}(x,
t)&=&\sum_{s=0}^{[\beta]}\tilde b^{(s)}\left(\frac{k}{L}x
+\frac{l}{L}2\pi j\right)\left(\frac l L\right)^s \frac{\left(t^s-(x-2\pi
j)^s\right)}{s!}.
\end{eqnarray*}
The last term may be written as the product of $x-2\pi j-t$ with a
  polynomial of degree less than $k$ in $t$. Because of the
  condition on the moments of $f$, the corresponding term is zero.

Let us then consider the second term, which is bounded by $C
|x-2\pi j|^\beta \Vert b\Vert_{\Lambda_\beta(\T)}$, where $C$ is an admissible
constant. To get a bound for the corresponding term in $\overline
T(f)(x)$, we use the fact that the Hilbert transform of $f$, at
the point $x-2\pi j$, can be bounded by
$Cr^{[\beta]+1-\beta}/|j|^{[\beta]+2}$. We conclude, using the
fact that $\sum_{|j|\geq 1}|j|^{\beta-[\beta]-2}$ is finite, and
$r^{[\beta]+1-\beta}$ is bounded.

Let us finally consider the first term, which  is bounded by $C
r^\beta\Vert b\Vert_{\Lambda_\beta(\T)}$, with $C$ an admissible constant. We
proceed exactly as in the proof of (\ref{bar-1}), using the fact that
the $L^1$-norm of $f$ is bounded by $r^{-\beta}$. We may as
well replace ${t +2\pi j-x}$ by $2\pi j$ in the last integral in
the denominator. Moreover, one  easily verifies that $e^{i\frac
\mu l x}\tilde B_j^{(1)}(x, t)$ can be written as $e^{i\frac \mu L
x}e^{-2\pi i\frac {\mu j} L }c\left(\frac kL (x+2\pi \frac lk
j),t\right)$, where the function $c$ is periodic in the first
variable $x$, and is bounded by $C r^\beta\Vert b\Vert_{\Lambda_\beta(\T)}$,
with $C$ an admissible constant which is independent of $t$. After
this point, the proof is exactly the same as the one of
(\ref{bar-1}), considering the sum $$\sum_{1\leq |j|<L/2}
c\left(\frac kL \left(x+2\pi \frac lk j\right),t\right)\frac {e^{-2\pi i\frac
{\mu j} L}}j$$  We will not repeat the proof here, and leave the
details for the reader. This finishes the proof of Proposition
\ref{bar}, as well as the proof of Theorem \ref{main}.
\end{proof}

\end{document}